\documentclass[twocolumn, twoside, 10pt, cefc]{IEEEtran}

\usepackage{graphicx,color,epsf,bm}
\usepackage{amsmath,amsfonts,amssymb,amscd,bm}
\usepackage{flushend}

\usepackage[colorlinks=true, urlcolor=blue]{hyperref}

\begin{document}
\title{\fontsize{17}{20}\bf
Inter-winding Distributed Capacitance and \\Guitar Pickup Transient Response}%
\author{
\IEEEauthorblockN{P. Robert Kotiuga\IEEEauthorrefmark{1}}
		\\
		\vspace{1ex}
    {\IEEEauthorrefmark{1}\small Boston University, Dept of ECE, 8 Saint Mary's Street, Boston MA 02215 USA}
   
}
\markboth{$>$ 5(\MakeLowercase{c, g,h, i}), 6(\MakeLowercase{e}) or 11(\MakeLowercase{m})$<$}{CEFC 2014 digest}
\maketitle

\begin{abstract}
Simple RLC circuit models of guitar pickups do not account for audible features that characterize the pickup. Psycho-acoustic experiments reveal that any acoustically accurate model has to reproduce the first 30 milli-seconds of the transient response with extreme precision. The proposed model is impractical for simple-minded model reduction or brute force numerical simulations yet, by focusing on modeling electromagnetic details and exposing a connection to spectral graph theory, a framework for finding the transient response to sufficient accuracy is exposed.
\end{abstract}
\begin{IEEEkeywords}Distributed capacitance, spectral graph theory, transient response.\end{IEEEkeywords}

\section{Introduction}

The pickup model proposed in this paper has a voltage and a current associated with each turn. Although this involves thousands of variables, the goal is to devise a model reduction strategy which accounts for the audible aspects of the transient response which characterize the ``tone" \cite{Hunter} of the pickup.  The ``cut-bell" psycho-acoustic experiments of Pierre Schaeffer reveal that any acoustically accurate model has to reproduce the first 30 milli-seconds of the transient response with extreme precision.\footnote{Although the brain cannot extract much information based on less than three milli-seconds of transient response due to the bandlimited nature of the hearing process, it processes transient information quite thoroughly based on less than 50 milliseconds of the waveform's ``attack".  Chapter 12, ``Anamorphoses Temporelles I: Timbres et Dynamiques", of Schaeffer's treatise\cite{Schaeffer}, summarizes this work. There are pointers to earlier observations of F. Winckel in this regard, but Schaeffer seems to have provided the key experimental confirmation and the readily reproducible methodology.}

The proposed pickup model is not straightforward, but has precedent in the analysis of very fast transients in multiwinding power transformers, a field with a history of detailed mathematical modeling that goes back a century \cite{Bjerkan} and where the literature on electromagnetics-based state-space analysis of transient oscillations, backed by rigorous numerical analysis, goes back at least 40 years \cite{Degeneff}. Problems arising in the context of power transformers also shed light on considerations of ``potting pickups" to suppress microphonics, and of disregarding eddy currents on ferromagnetic parts. Unfortunately, this analogy does not resolve the issues involved in modeling ``scatterwound" pickups since the model reduction techniques made in the context of power transformers oversimplify the pickup problem. For instance, for the purpose of capacitance modeling, the turns of transformer winding can be modeled as sheet currents \cite{FerHenr, Hosseinietal, Gomezetal}.  Also, unlike power transformers, pickup wire is typically 42 AWG, so skin and proximity effects are non-issues well into the MHz range. The key common aspect is the neccessity of starting with an accurate multiconductor transmission line model \cite{Paul}, independent of questions pertaining to model reduction. 

To fix ideas, consider three questions which one might hope to answer in the context of analysing a pickup:
\begin{enumerate}
\item Deriving an accurate transfer function relating the motion of the string to the output voltage of the pickup;
\item Finding the input impedance of the pickup.
\item Find the (short-circuit) internal resonances of the pickup.
\end{enumerate}
The third problem is a stepping stone to answering the first two questions. In this paper we do not attempt to answer the first question in full generality. This is because an accurate model of pickup would involve, for instance, the eddy-currents on the pole pieces of the pickup. Such eddy currents lead to ``wolf sound". For a properly adjusted guitar, string heights above pole pieces are adjusted to make it imperceptible and so we do not need to model the eddy-currents. Similarly, we do not dwell on the second problem since the input impedance  only matters in the context of ensuring that the impedance of the guitar and cable, as seen from the follow-on stage of amplification, is negligible. Hence, in this paper we focus on the role of the of distributed capacitance in the context of the third problem. This in turn will lead us to spectral graph theory \cite{Brualdi}, \cite{Chung}, a formalization of the well-studied link between nodal analysis, the graph Laplacian and the spectrum of the adjacency matrix of a graph.  

\section{A Sketch of an Idea}

This paper exposes the connection between the natural frequencies of the pickup and the eigenvalues of the adjacency matrix. This connection can be sketched as follows. Consider the following linear system arising from the discretization of multi-conductor transmission lines which are wrapped around to form the turns of a winding. It involves square matrices of order $n$, where $n$ is the number of turns:
\begin {equation}\label{prototype}
(\mathsf{\Gamma^{2}}-\mathsf{Circ})\mathsf{x}=0.
\end {equation}
Here, $\mathsf{Circ}$ is a circulant matrix and $\mathsf{x}$ are state variables. That is, a matrix whose (k+1)-st row is the k-th row shifted one entry right and wrapped around. In general any circulant matrix is a linear combination of the first $n$ powers of the shift matrix $\mathsf{S}$. In other words, $\mathsf{Circ}$ is a polynomial in $\mathsf{S}$.
Following the theory of the DFT, we know $\mathsf{S}$ is diagonalized by the unitary matrix $\mathsf{V}$, a Vandermonde matrix of roots of unity, normalized by $n^{{-1}/{2}}$,  and the eigenvalues of $\mathsf{S}$ are the $n^{th}$ roots of unity. That is, $\mathsf{S}=\mathsf{V}^{\dagger}\Lambda_{\mathsf{S}}\mathsf{V}$, where $\Lambda_{\mathsf{S}}$ is the diagonal matrix of roots of unity, and $\dagger$ denotes the Hermitian conjugate. Furthermore, all powers of $\mathsf{S}$ are also diagonalized by $\mathsf{V}$, as are polynomials and analytic functions of $\mathsf{S}$. So, 
\begin {equation}\label{circ}
\mathsf{Circ}=P_{1}(\mathsf{S})=\mathsf{V}^{\dagger}P_{1}(\Lambda_{\mathsf{S}})\mathsf{V}
\end {equation} 
The square of the propagation matrix $\mathsf{\Gamma^{2}}$ is modeled as a polynomial in the adjacency matrix $\mathsf{A}$ of a graph associated to the manner in which the turns of the pickup winding were wound, with the angular frequency $\omega$ as a parameter:
\begin {equation}\label{propagation}
\mathsf{\Gamma^{2}}=P_{2}(\mathsf{A},\omega)
\end {equation}
Eq. (\ref{propagation}) builds on the relationship between the ``Laplacian of a graph" where the degree of any node is independent of the node index, and the nodal analysis matrix of a network having the same underlying graph. In this paper the node degree is dependent on the node index at the boundary of the winding and an additional diagonal matrix will be introduced to account for the discrepancy. Substitute (\ref{circ}) and (\ref{propagation}) into (\ref{prototype}) to get:
 \begin {equation}\label{sketchy}
(P_{2}(\mathsf{A},\omega)-P_{1}(\mathsf{S}))\mathsf{x}=0 \;\;\;\;\;\; \text{or},
\end {equation} 
\begin {equation}\label{sketchier}
(P_{2}(\mathsf{A},\omega)-\mathsf{V}^{-1}P_{1}(\Lambda_{\mathsf{S}})\mathsf{V})\mathsf{x}=0
\end {equation} 
We would expect to find nontrivial solutions by setting the determinant to zero to obtain a polynomial equation in $\omega$ whose roots give the natural frequencies. Unfortunately, for the model at hand all of the matrices are singular with a common null-space and common null vector. Hence the matrix is singular for all values of $\omega$. In section six, we will obtain the correct  equation by using the decomposition in (\ref{circ}) to project onto the orthogonal complement of the null-vector, and then form the determinantal equation for the natural frequencies.

\section{The Network Model}

Consider a single-coil magnetic pickup on an electric guitar, placed in 3-D Euclidean space with coordinates ($x,y,z$) such that the $z$-axis is perpendicular to the fretboard. The pickup in question is centered on the $z$-axis with its center at $z$=0. From a field theoretic vantage, the subtle aspect of the guitar pickup model is that each turn of the winding links a nontrivial time-varying magnetic field, yet one describes the electric field exterior to the conducting coil in terms of an electric scalar potential. This is accomplished within the convex hull of the pickup by assuming that the time rate of change of magnetic flux perpendicular to any plane containing the $z$-axis is negligible. This assumption enables one to define an electric scalar potential for the components of the electric field lying in any plane containing the $z$-axis\footnote{ In principle, the planarity condition in this assumption is only valid for a circular coil. A more precise formulation of this assumption would involve a singular foliation of the convex hull of the pickup with the all of the leaves bordering on the z-axis. A precise formulation of the assumption is not required since only the existence of the scalar potential is required to formulate the network equations.}. However, the time-varying magnetic flux linking the pickup coil prevents this electric scalar potential from describing the “circumferential” components. To describe these components of the electric field within the conductor, we integrate the electric field along the conductor and evaluate the value of the line integral every time one passes through a given half-plane whose boundary is the z-axis. Thus let
\begin {equation}
\int_{c_i}\boldsymbol{E}\cdot \boldsymbol{dl}= V_i -V_0
\end {equation}
where $c_i$ comprises the the first i turns of the winding. The voltage drop around the $i$-th turn is then given by $V_i -V_{i-1}$. We can compile these voltages into a vector $\boldsymbol{V}$. Similarly, if $I_i$ is the current entering the $i$-th winding, these currents can be compiled into a vector $\boldsymbol{I}$. Let $n$ be the number of turns. Taking the index modulo $n$, the equations for the short circuited pickup are given by 
\begin {equation}\label{KVL}
(\mathsf{S}-\mathsf{I})\boldsymbol{V}= -(\mathsf{R}+j\omega\mathsf{L})\boldsymbol{I} \;\;\;\;\;\; \text{and},
\end {equation}
\begin {equation}\label{KCL}
(\mathsf{S}-\mathsf{I})\boldsymbol{I}= -j\omega\mathsf{C}\boldsymbol{V},
\end {equation}
where $\mathsf{I}$ is the identity matrix,  $\mathsf{S}$ is the shift matrix, and the remaining impedance and admittance matrices are as one would expect if one took a length of multiconductor transmission line \cite{Paul} and connected the end of the $k$-th turn to the beginning of the $(k+1)$-st. To preserve the cyclic symmetry, we avoided grounding any node and so the capacitance matrix is necessarily singular with null vector $\mathsf{1}$, where $\mathsf{1}$ is the vector all of whose entries are $1$. Furthermore, the matrix $(\mathsf{S}-\mathsf{I})$ has the same null vector. This approach is formalized in the method of the indefinite impedance matrix\footnote{ A textbook exposition of the indefinite impedance matrix method and its relationship to cutset analysis can be found in section 3.7 of Balabanian and Bickart's text \cite{BalanBick}. Note that for a passive network this impedance matrix has one zero eigenvalue for each connected component of the network and the remaining eigenvalues have positive real parts}.

\section{An Analog of the Telegraphist's Equation via Simple Models of $\mathsf{R}$ and $\mathsf{L}$ Matrices}

If we multiply (\ref{KVL}) by $j\omega\mathsf{C}$ we obtain:
\begin {equation}\label{something}
j\omega\mathsf{C}(\mathsf{S}-\mathsf{I})\boldsymbol{V}= -j\omega\mathsf{C}(\mathsf{R}+j\omega\mathsf{L})\boldsymbol{I}
\end {equation}
 Add and subtract $j\omega(\mathsf{S}-\mathsf{I})\mathsf{C}\boldsymbol{V}$ to the l.h.s. of (\ref{something}) and simplify, to get a telegraphist-like equation, plus a commutator term:
\begin {equation}\label{1sttelegraphist}
(-j\omega\mathsf{C}(\mathsf{R}+j\omega\mathsf{L})+ (\mathsf{S}-\mathsf{I})^2)\boldsymbol{I}=j\omega(\mathsf{C}\mathsf{S}-\mathsf{S}\mathsf{C})\boldsymbol{V}
\end {equation}
Had the capacitance matrix been invertible, and if  $\mathsf{C}$ and $(\mathsf{S}-\mathsf{I})$ commuted, we could hope to eliminate $\boldsymbol{V}$ using (\ref{KCL}).  We will see that if we formulate our equations in terms of voltages, the current can be eliminated yielding an analog of the telegraphist equation, thus avoiding these difficulties. Additional structure in the $\mathsf{R}$, $\mathsf{L}$ and $\mathsf{C}$ matrices leads to useful commutation relations. This is easy to verify for the  $\mathsf{R}$ and $\mathsf{L}$ matrices. The commutation relations will ensure that the analog of the telegraphist's equations, when formulated in terms of voltages, are free of undesirable commutators. 

The coil in a guitar pickup consists of 5,000 to 8,000 turns of copper wire (typically 42 AWG) with a typical DC resistance of 7.5 k$\Omega$.
The matrix $\mathsf{R}$ is modeled as a multiple of the identity matrix and is the simplest to model. If $R_\text{DC}$ is the DC resistance of the pickup, and the resistance, $R_0$ of one turn is assumed to be the average, independent of the turn, then $R_0$ is $R_{DC}$ divided by $n$. Hence,
\begin {equation}\label{R}
\mathsf{R}=R_{0}\mathsf{I} \;\;\text{and } \text{Trace}(\mathsf{R})=R_{DC}\;\;\text{where } R_{0}= \frac{R_{DC}}{n},
\end {equation}
The inner turns of the winding are shorter than the outer ones, but this is ignored in the analysis. 
Since $\mathsf{R}$ is a multiple of the identity matrix, it commutes with any compatible matrix.

 Let $M$ be the mutual inductance between two turns and assume this value is independent of the index of the turn. Let ${\alpha}$ be a nonnegative constant expressing the deviation from perfect coupling. The inductance matrix model we use is
\begin {equation}\label{L}
\mathsf{L}=M({\alpha}\mathsf{I}+(\mathsf{1}\otimes \mathsf{1})).
\end {equation}
Here, $\mathsf{1}\otimes \mathsf{1}$ the outer product of the vector $\mathsf{1}$ with itself and again, $\mathsf{1}$ is the vector all of whose entries are $1$. Since the capacitance matrix is an admittance matrix arising from nodal analysis, and a ground node has not been selected, it is symmetric positive semi-definite with a one-dimensional null space. This is reflected in the fact that the row sums or column sums of any row or column add up to zero. The null vector is precisely the vector $\mathsf{1}$ which occurs in the outer product in the inductance matrix. In other words, $\mathsf{C}\mathsf{1}=\mathsf{0}$, so 
\begin {equation}\label{Cnull}
\mathsf{C}(\mathsf{1}\otimes \mathsf{1})=\mathsf{0},
\end {equation}
and, by (\ref{L}) we have
\begin {equation}\label{CL}
\mathsf{C}\mathsf{L}=\mathsf{C}(M{\alpha}\mathsf{I}+(\mathsf{1}\otimes \mathsf{1}))={\alpha}M\mathsf{C}.
\end {equation}
Furthermore, since both $\mathsf{L}$ and $\mathsf{C}$ are symmetric the same result holds for the commuted product, so that
\begin {equation}\label{CLLC}
\mathsf{C}\mathsf{L}={\alpha}M\mathsf{C}=\mathsf{L}\mathsf{C}. 
\end {equation}
This argument shows that $\mathsf{L}$ commutes with any symmetric matrix whose null-space includes the vector $\mathsf{1}$ and easily extends to nonsymmetric matrices where $\mathsf{1}$ is included in both the nullspaces of the matrix and its transpose. One example of such a matrix is $(\mathsf{S}-\mathsf{I})$, and one can easily verify that
\begin {equation}\label{SLLS}
\mathsf{(\mathsf{S}-\mathsf{I})}\mathsf{L}={\alpha}M\mathsf{(\mathsf{S}-\mathsf{I})}=\mathsf{L}\mathsf{(\mathsf{S}-\mathsf{I})}. 
\end {equation}
From eqs (\ref{KVL}, \ref{KCL}) and commutation relations (\ref{CLLC}, \ref{SLLS}) we obtain:
\begin {eqnarray*}
(\mathsf{S}-\mathsf{I})^2\boldsymbol{V} & = & -(\mathsf{S}-\mathsf{I})(\mathsf{R}+j\omega\mathsf{L})\boldsymbol{I}\\
		& = & -(\mathsf{S}-\mathsf{I})(\mathsf{R}+j\omega M({\alpha}\mathsf{I}+(\mathsf{1}\otimes \mathsf{1})))\boldsymbol{I}\\
		& = & -(R_{0}+j\omega M{\alpha})(\mathsf{S}-\mathsf{I})\boldsymbol{I}\\
		& = & j\omega\mathsf{C}(R_{0}+j\omega M{\alpha})\boldsymbol{V}, \text{or}\\ 	
\end {eqnarray*}
\begin {equation}\label{2ndtelegraphist}
((\mathsf{S}-\mathsf{I})^2 -j\omega\mathsf{C}(R_{0}+j\omega M{\alpha}))\boldsymbol{V}= 0.
\end {equation}
This is a commutator-free analog of the telegraphist equation. 

\section{Winding Induces an Ordering on Turns; Adjacency and Distributed Capacitance Matrices}

In order to model the inter-winding capacitance matrix we introduce a directed graph, $G_d$, as follows. Assume that the turns are ordered by the order in which they were wound so that the $i$-th index in the matrix equation corresponds to the $i$-th turn wound. This is our definition of a ``preferred ordering''. The edges of $G_d$ correspond to neighboring turns, inducing a ``preferred orientation'' on the edges of $G_d$ dictated by the order the turns were wound. It follows that the edges of $G_d$ are directed from higher to lower node index.
Associated to $G_d$ is the undirected graph, $G_u$, obtained by ignoring the orientations of the edges. We associate the adjacency matrix $\mathsf{A}$ to $G_u$ by

\begin{equation}\label{A}
        (\mathsf{A})_{i,j} = \left\{\begin{array}{cl}
        0 & \text{if $i=j$ or no edge connects nodes $i$ and $j$},\\
        1 & \text{otherwise}.
        \end{array}\right. 
\end{equation}

If we assume that, locally, the cross-section of the winding appears to be a hexagonal closest packing, then any row or column corresponding to a non-boundary winding has exactly six nonzero off-diagonal entries. The existence of $G_d$ ensures that, on average, every node of $G_d$ corresponding to a non-bounday turn has three incoming and three outgoing edges. 

From our discussion of the winding's cross-section, we have the following model of the inter-winding capacitiance matrix,
\begin {equation}\label{C}
\mathsf{C}=C_{0}(d\mathsf{I}-\mathsf{A}+\mathsf{B})\;\text{and }  \text{Trace}(\mathsf{C})=(n_\text{int}d+ n_\text{b}d_\text{bav})C_{0},
\end {equation}
since $ \text{Trace}(\mathsf{A})=0$. Here, $C_{0}$ is the capacitance one has between two turns of enameled  copper wire wound beside each other. This parameter is easily deduced by conformal mapping once the diameter of the enamel coating is known. The matrix $\mathsf{B}$ is a diagonal matrix with a relatively small percentage of nonzero diagonal entries corresponding to boundary turns. The value of the nonzero entry corresponding to a given boundary node is the degree of the node minus six. In (\ref{C}), $n_{int}$ and $n_{b}$ are respectively the number of internal and boundary nodes, $d$ is the degree of the internal nodes, assumed constant (it is six but we leave it as a symbol to suggest that $G_u$ is a regular graph), and $d_{bav}$ is the average degree of a boundary node. 

Although the incidence data of $G_d$ appears in the statement of Kirchhoff's laws, it is only $G_u$ which appears explicitly in the nodal equations through the appearance of $\mathsf{A}$. However, in the case at hand, the orientation information in the graph $G_d$ is implicit in the nodal equations once the preferred ordering is introduced. This is because the $k$-th edge originates on a node whose index $j$ is less or equal to $k$ and terminates on a node of index lower than $j$. In this way, the orientation information of $G_d$ can then be reconstructed from the sparsity pattern of $\mathsf{A}$, and it is in this manner that winding pattern affects the eigenvalue distribution of $\mathsf{C}$ through $\mathsf{A}$.

We have seen that $\mathsf{R}$, $\mathsf{L}$ and $\mathsf{C}$ are a set of commuting matrices and that $\mathsf{R}$, $\mathsf{L}$ and $(\mathsf{S}-\mathsf{I})$ are another such set. To see why $\mathsf{C}$ and $(\mathsf{S}-\mathsf{I})$ are not, use (\ref{C}) to reduce the commutator on the r.h.s. of (\ref{1sttelegraphist}) to those between $\mathsf{A}$, $\mathsf{B}$ and $\mathsf{S}$: 
\begin {equation}\label{CSCommutator}
\mathsf{C}\mathsf{S}-\mathsf{S}\mathsf{C}={C_{0}}((\mathsf{B}-\mathsf{A})\mathsf{S}-\mathsf{S}(\mathsf{B}-\mathsf{A})).
\end {equation}
Finally, substituting (\ref{R},\ref{L},\ref{CL}) into (\ref{1sttelegraphist}), reduces it to: 
\begin {equation}\label{homo nodal eq}
((\mathsf{S}-\mathsf{I})^2-j\omega\mathsf{C}(R_{0}+j\omega{\alpha}M))\boldsymbol{I}=j\omega(\mathsf{C}\mathsf{S}-\mathsf{S}\mathsf{C})\boldsymbol{V}.
\end {equation}
Eq. (\ref{CSCommutator}) reveals how formulating our telegraphist equations in terms of $\boldsymbol{I}$ leads to unavoidable commutators. And so we return to (\ref{2ndtelegraphist}) in order to relate the natural frequencies of this model to the $\mathsf{C}$ matrix via the sparsity structure of the adjacency matrix. If we substitute (\ref{C}) into (\ref{2ndtelegraphist}) we get
\begin {equation}\label{2ndteleWC}
((\mathsf{S}-\mathsf{I})^2 -j\omega C_{0}(R_{0}+j\omega M{\alpha})(d\mathsf{I}-\mathsf{A}+\mathsf{B}))\boldsymbol{V}= 0
\end {equation}
or
\begin {equation}\label{2ndteleWgamma}
((\mathsf{S}-\mathsf{I})^2 -\gamma^{2}(d\mathsf{I}-\mathsf{A}+\mathsf{B}))\boldsymbol{V}= 0;
\end {equation}
\begin {equation}\label{gamma}
\gamma^{2}=j\omega C_{0}(R_{0}+j\omega M{\alpha})
\end {equation}

\section{The tie to spectral graph theory}

Eq. (\ref{2ndteleWgamma}) cannot be a conventional eigenvalue problem for $\gamma^{2}$ since it turns out that the underlying matrices are singular for all values of $\gamma^{2}$. To resolve this problem, we note that the vector $\mathsf{1}$  has appeared in many contexts. Specifically, it is the null vector of $\mathsf{C}=C_{0}(d\mathsf{I}-\mathsf{A}+\mathsf{B})$, the null vector of $(\mathsf{S}-\mathsf{I})$ and the eigenvector of $\mathsf{S}$ corresponding to the eigenvalue $1$. One could toss out a row and column from the system, as when one grounds a node in nodal analysis, but a  better approach is to project the system into the orthogonal complement of the vector $\mathsf{1}$. Recall that the Vandermonde matrix of $n$-th roots of unity diagonalizes any circulant matrix and that the eigenvalues of the shift matrix are just the roots of unity. So,  let $V^{\dagger}$ denote the Hermitian conjugate of $V$ and let 
\begin {equation}\label{VanderMon_n}
\Lambda_n = \text{diag}( [ \Omega_{n}^{k}]_{k=o}^{n-1} )\;\;\text{where} \;\Omega_{n}^{k}= e^{\frac{2\pi ki}{n}} . 
\end {equation}
The first term in (\ref{2ndteleWgamma}) is a circulant matrix which is diagonalized by the similarity transfomation of section two:
\begin {equation}\label{SminusIsquared}
\mathsf{(\mathsf{S}-\mathsf{I})^2}=\mathsf{V}^{\dagger}((\Lambda_{n}-\mathsf{I})^2)\mathsf{V}= \mathsf{V}^{\dagger}((\widehat{\Lambda_{n}})\mathsf{V},
\end {equation}
where by (\ref{VanderMon_n}) we explicitly have
 \begin {equation}\label{LambdaNhat}
\widehat{\Lambda_{n}}= \text{diag}( [e^{\frac{2\pi ki}{n}}(cos(\frac{2\pi ki}{n})-1) ]_{k=o}^{n-1} .
\end {equation}
Substituting (\ref{SminusIsquared}) into (\ref{2ndteleWgamma}), we can write
\begin {equation}\label{2ndteleWgammaNV}
(\mathsf{V}^{\dagger}\widehat{\Lambda_{n}}\mathsf{V} -\gamma^{2}(d\mathsf{I}-\mathsf{A}+\mathsf{B}))\boldsymbol{V}= 0.
\end {equation}
Although this is still singular for all values of $\gamma^{2}$, the similarity transformation induced by $\mathsf{V}$ points the way to projecting into the orthogonal complement to the vector $\mathsf{1}$. Rewrite (\ref{2ndteleWgammaNV}) as
\begin {equation}\label{2ndteleWgammaNvMassaged}
(\widehat{\Lambda_{n}}-\gamma^{2}(d\mathsf{I}-\mathsf{V}(\mathsf{A}-\mathsf{B})\mathsf{V}^{\dagger}))\mathsf{V}\boldsymbol{V}= 0.
\end {equation}
There is one row and one column which are explicitly zero and so we can now explicitly perform the projection. Let
 \begin {equation}\label{LambdaNminus1tilde}
\widetilde{\Lambda_{n-1}}= \text{diag}( [e^{\frac{2\pi ki}{n}}(cos(\frac{2\pi ki}{n})-1) ]_{k=1}^{n-1}. 
\end {equation}
That is, with $k\ne 0,n$, we restricted ourselves to the nonzero eigenvalues of $(\mathsf{S}-\mathsf{I})$. Let $\mathsf{V}^{\dagger}_{Pr}$ be the $n\times(n-1)$ matrix consisting of the eigenvectors of $(\mathsf{S}-\mathsf{I})$ corresponding to nonzero eigenvalues. Then, 
\begin {equation}\label{unit}
\mathsf{V}^{\dagger}\mathsf{V}=\mathsf{V}\mathsf{V}^{\dagger}=\mathsf{I}_{n{\times}n} \;\;\text{but }
\end {equation}
\begin {equation}\label{proj}
\mathsf{V}_{Pr}\mathsf{V}_{Pr}^{\dagger}=\mathsf{I}_{(n-1){\times}(n-1)} \; \text{and} \; \mathsf{V}_{Pr}^{\dagger}\mathsf{V}_{Pr}={Proj}\perp{\mathsf{1}}
\end {equation}
The determinantal equation is now expressed as:
\begin {equation}\label{det2ndteleWgammaNvMassaged}
\text{det}(\widetilde{\Lambda_{n-1}}-\gamma^{2}(d\mathsf{I}_{(n-1){\times}(n-1)}-\mathsf{V}_{Pr}(\mathsf{A}-\mathsf{B})\mathsf{V}_{Pr}^{\dagger}))= 0
\end {equation}
This is the definitive equation for subsequent developments. However, given the definition of the capacitance matrix, this equation could have been written in the more intuitive form
\begin {equation}\label{det2ndteleWgammaNvIntuitive}
\text{det}(\widetilde{\Lambda_{n-1}}+\frac{\gamma^{2}}{C_{0}}\mathsf{V}_{Pr}\mathsf{C}\mathsf{V}_{Pr}^{\dagger})= 0.
\end {equation}

Hence, (\ref{det2ndteleWgammaNvMassaged}) reduces the problem to a modified eigenvalue problem for a (diagonal) perturbation of the (vanishing diagonal) adjacency matrix. By (\ref{gamma}) any information about the spectrum translates into information about the collection of natural frquencies since they are related to the eigenvalues through a quadratic equation.  The adjacency matrix eigenvalue problem is the focus of ``spectral graph theory'', and there are many ways to exploit this connection. This provides precise estimates for the transient response and a framework for model reduction in the presence of eigenvalue clusters.

\section{Summary and Conclusions}

In this paper we shed light on elusive aspects of electromagnetic modeling of  electric guitar pickups. Specifically, the ability of an expert to infer the manner in which the pickup was wound by listening.  Psycho-acoustic experiments reveal that acoustically accurate models have to reproduce the first 30 milliseconds of the transient response with extreme precision since the brain makes inferences about the pickup's tone within this timeframe. Conventional wisdom dictates that one could take the natural frequencies which are smallest in absolute value and perform a model reduction. This however is tricky when the natural fequencies cluster. The problem is difficult because the winding pattern of the pickup affects the clustering of natural frequencies in a subtle manner which is audibly discernable as the tone of the pickup, yet gets lost in a simple-minded model reduction scheme. 

By exposing a connection to spectral graph theory, (\ref{2ndteleWgamma},\ref{gamma},\ref{det2ndteleWgammaNvMassaged}) yield a framework for analyzing the ``attack" of the transient response to sufficient detail to reveal the role of the ordering of the turns in the winding. This model is insensitive to how the pickup was wound with the exception of the inter-winding capacitance matrix which describes the capacitance between the turns of the winding.

\end{document}